\documentclass[12pt]{amsart}  
\usepackage{amssymb,amsthm}  
\usepackage{amscd}

\newcommand{\Z}{\mathbb{Z}}  
\newcommand{\N}{\mathbb{N}}  
  
\newcommand{\ab}{{\bf a}}

\newcommand{\xb}{{\bf x}}

\DeclareMathOperator{\chara}{char}
  
\DeclareMathOperator{\codim}{codim}

\DeclareMathOperator{\Gin}{Gin}

\DeclareMathOperator{\Ker}{Ker}

\DeclareMathOperator{\projdim}{proj\,dim}

\DeclareMathOperator{\reg}{reg}

\DeclareMathOperator{\supp}{supp}

\DeclareMathOperator{\Tor}{Tor}

\DeclareMathOperator{\tensor}{\otimes}

\DeclareMathOperator{\pnt}{\raise 0.5mm \hbox{\large\bf.}}  
\DeclareMathOperator{\lpnt}{\hbox{\large\bf.}}

\newtheorem{Thm}{\bf Theorem}[section]   
\newtheorem{Lem}[Thm]{\bf Lemma}  
\newtheorem{Cor}[Thm]{\bf Corollary}  
  
\newtheorem{Rem}[Thm]{\bf Remark}

\newtheorem{Conj}[Thm]{\bf Conjecture}

\title{Note on bounds for multiplicities}  

\author{Tim R\"omer} 
\address{FB Mathematik/Informatik, Universit\"at Osnabr\"uck, 49069 Osnabr\"uck, Germany}
\email{troemer@mathematik.uni-osnabrueck.de}

\begin{document}  
  
\begin{abstract}  
Let $S=K[x_1,\ldots,x_n]$ be a polynomial ring 
and $R=S/I$ be a graded $K$-algebra where $I \subset S$ is a graded ideal.
Herzog, Huneke and Srinivasan 
have conjectured that the multiplicity of $R$ is bounded
above by a function of the maximal shifts in 
the minimal graded free resolution of $R$ over $S$.
We prove the conjecture in the case that $\codim(R)=2$ which
generalizes results in \cite{Go} and \cite{HeSr98}.
We also give a proof for the bound 
in the case in which $I$ is componentwise linear. 
For example, stable and squarefree stable ideals belong 
to this class of ideals.
\end{abstract}
  
\maketitle 

%
%
%
\section{Introduction}  
Let $S=K[x_1,\ldots,x_n]$ be the polynomial ring with $n$ variables over a field $K$  
equipped with the standard grading 
by setting $\deg(x_i)=1$. 
Let $I\subset S$ be a graded ideal and $R=S/I$ be a standard graded $K$-algebra.
Consider the minimal graded free resolution of $R$:
$$
0
\to 
\bigoplus_{j\in \Z} S(-j)^{\beta^S_{p,j}(R)} 
\to
\cdots
\to
\bigoplus_{j\in \Z} S(-j)^{\beta^S_{1,j}(R)}
\to S \to 0 
$$
where we denote with  
$\beta^S_{i,j}(R)=\dim_K \Tor^S_i(R,K)_{j}$ the graded Betti numbers of $R$
and $p=\projdim(R)$ is the projective dimension of $R$.
The ring $R$ is said to have a {\it pure} resolution if
at each step there is only a single degree, i.e. 
the resolution has the following shape:
$$
0
\to 
S(-d_p)^{\beta^S_{p}(R)} 
\to
\cdots
\to
S(-d_1)^{\beta^S_{1}(R)} 
\to S \to 0 
$$
for some numbers $d_1,\ldots,d_p$.

Let $e(R)$ denote the multiplicity of $R$. 
Huneke and Miller proved in \cite{HuMi85} the following formula:

\begin{Thm}
Let $R$ be a Cohen-Macaulay ring with a pure resolution. Then
$$
e(R) = (\prod_{i=1}^{p} d_i)/p!.
$$
\end{Thm}

More general, for $1\leq i \leq \projdim(R)$ we define
$$
M_i=\max\{j \in \Z \colon \beta^S_{i,j}(R)\neq 0\}
\text{ and }
m_i=\min\{j \in \Z \colon \beta^S_{i,j}(R)\neq 0\}.
$$
The multiplicity conjecture of 
Huneke and Srinivasan 
is:

\begin{Conj}
\label{conj1}
Let $I \subset S$ be a graded ideal, $R=S/I$ be Cohen-Macaulay and $p=\projdim(R)$. Then
$$
(\prod_{i=1}^{p} m_i)/p!\leq e(R) \leq (\prod_{i=1}^{p} M_i)/p!.
$$
\end{Conj}

Herzog and Srinivasan proved this conjecture in \cite{HeSr98}
for several types of ideals:
complete intersections, 
perfect ideals with quasipure resolutions (i.e. $m_i(R)\geq M_{i-1}(R)$ for all $i$), 
perfect ideals of codimension 2, 
codimension 3 Gorenstein ideals generated by 5 elements 
(the upper bound holds for all codimension 3 Gorenstein ideals),
codimension 3 Gorenstein monomial ideals with at least one generator of smallest possible degree (relative to the number of generators),
perfect stable ideals (in the sense of Eliahou and Kervaire \cite{ElKe}),
perfect squarefree strongly stable ideals (in the sense of Aramova, Herzog and Hibi \cite{ArHeHi98}).
See also \cite{HeSr02} for related results. 
The lower bound fails to hold in general if $R$ is not Cohen-Macaulay 
(see \cite{HeSr98} for a detailed discussion). 
Herzog and Srinivasan conjectured in this case the following 
inequality:

\begin{Conj}
\label{conj2}
Let $I \subset S$ be a graded ideal, $R=S/I$ and $c=\codim(R)$. Then
$$
e(R) \leq (\prod_{i=1}^{c} M_i)/c!.
$$
\end{Conj}

Since the codimension of a graded $K$-algebra is less or equal to the projective
dimension and for all $i$ we have that $M_i\geq i$, 
the inequality in Conjecture \ref{conj2} 
is stronger than the one of Conjecture \ref{conj1}.

Herzog and Srinivasan proved this conjecture in the cases
of stable ideals, squarefree strongly stable ideals and
ideals with a $d$-linear resolution, i.e.\ $\beta^S_{i,i+j}(R)=0$ for $j\neq d$.
Furthermore Gold \cite{Go} 
established Conjecture \ref{conj2} in the case of codimension
$2$ lattice ideals. This conjecture is also known to be true for so-called 
$\ab$-stable ideals 
(see Section \ref{compsect} for the definition) by Gasharov, Hibi and Peeva \cite{GAHIPE}
which generalizes the stable and squarefree stable case.

In the first part of this paper 
we show that Conjecture \ref{conj2} is valid 
for codimension 2 ideals.
This generalizes the cases of perfect codimension 2 ideals 
of Herzog and Srinivasan and codimension 2 lattice ideals of Gold.

For $d\geq 0$ let $I_{\langle d \rangle} \subseteq I$ 
be the ideal which 
is generated by all elements of degree $d$ in $I$.
Recall from \cite{HeHi99} 
that an ideal $I\subset S$ is called 
{\it componentwise linear} if for all $d\geq 0$
the ideal $I_{\langle d \rangle}$ has a $d$-linear resolution.
We show that the upper
bound for the multiplicity holds 
for componentwise
linear ideals which generalizes some of the known cases
since for example stable and squarefree stable ideals are componentwise linear.
We prove that $\ab$-stable ideals are componentwise linear and can
also deduce the conjecture in this case.

The author is grateful to Prof.\ Bruns and Prof.\ Herzog for inspiring 
discussions on the subject of the paper.

%
%
\section{Codimension 2 case}
Let $I\subset S$ be a graded ideal and $R=S/I$. 
In this section we prove Conjecture \ref{conj2} 
in the case that $\codim(R)= 2$.

The codimension 2 case is known 
if $R$ is Cohen-Macaulay:

\begin{Thm}[Herzog-Srinivasan \cite{HeSr98}]
\label{codim2_cm}
Let $I\subset S$ be a graded ideal and $R=S/I$ Cohen-Macaulay with 
$\codim(R)=2$. 
Then 
$$
e(R) \leq (M_1\cdot M_2)/2.
$$
\end{Thm}

Following \cite{ARHE00} (or \cite{TR} under the name filter regular element)
we call an element 
$x \in R_{1}$ {\it almost regular} for $R$
if 
$$
(0:_R x)_{a} = 0 \text{ for } a \gg 0.
$$
A sequence $x_1,\ldots,x_t \in R_{1}$ is an {\it almost regular sequence}
if for all $i \in \{1,\dots,t\}$ the element 
$x_i$ is almost regular for $R/(x_1,\ldots,x_{i-1})R$.
It is well-known that, provided $|K|=\infty$, 
after a generic choice of coordinates we can achieve that 
a $K$-basis of $R_{1}$ is almost regular for $R$. (See \cite{ARHE00} and \cite{TR} for details.)

If $\dim_K(R)=n$ and since neither the Betti numbers nor the multiplicity of $R$
changes by enlarging the field, we always may assume that $x_1,\ldots,x_n\in R_1$
is an almost regular sequence for $R$  
to prove Conjecture \ref{conj2}. 
In the following we will not distinguish between an element $x\in S_1$
and the image in $R_1$.

We use almost regular elements to reduce the problem to dimension zero.
At first we have to recall some properties of almost regular elements.
\begin{Lem}
\label{help_dim}
Let $I\subset S$ be a graded ideal and $R=S/I$. 
Let $x\in R_1$ be almost regular for $R$.  
If $\dim(R)>0$, then $\dim(R/xR)=\dim(R)-1$.
\end{Lem}
\begin{proof}
We have the exact sequence
$$
0 
\to 
(0:_R x)(-1)
\to
R(-1)
\overset{x}{\to}
R
\to
R/xR
\to
0.
$$
Since $(0:_R x)$ has finite length
and $\dim(R)>0$
we conclude that
$\dim(R/xR)=\dim(R)-1$.
\end{proof}

\begin{Lem}
\label{help_multi}
Let $I\subset S$ be a graded ideal and $R=S/I$. 
Let $x\in R_1$ be almost regular for $R$. 
Then
\begin{enumerate}
\item
If $\dim(R)>1$, then $e(R)=e(R/xR)$.
\item
If $\dim(R)=1$, then $e(R)\leq e(R/xR)$.
\end{enumerate}
\end{Lem}
\begin{proof}
Again we have the exact sequence
$$
0 
\to 
(0:_R x)(-1)
\to
R(-1)
\overset{x}{\to}
R
\to
R/xR
\to
0.
$$
In the case $\dim(R)>1$ we get $e(R)=e(R/xR)$, because
$(0:_R x)$ has finite length.
If $\dim(R)=1$, then 
$$
e(R)=e(R/xR)-l((0:_R x))\leq e(R/xR).
$$
\end{proof}

Let 
$K_{\lpnt}(k;R)$ denote the Koszul complex 
and $H_{\lpnt}(k;R)$ denote the Koszul homology of $R$
with respect to $x_1,\ldots,x_k \in S$ (see \cite{BRHE98} for details).
Note that $K_{\lpnt}(k;R)=K_{\lpnt}(k;S)\tensor_S R$
where $K_{\lpnt}(k;S)$ is the exterior algebra on
$e_1,\ldots,e_k$ with $\deg(e_i)=1$ 
together with a differential $\partial$
induced by $\partial(e_i)=x_i$.
For a cycle $z \in K_{\lpnt}(k;R)$ we denote with $[z]\in H_{\lpnt}(k;R)$
the corresponding homology class. 
For $k=0$ we set $H_0(0;R)=R$. Then
there are long exact sequences relating the Koszul homology groups:
$$
\cdots \to H_i(k;R)(-1)\overset{x_{k+1}}{\to} H_i(k;R)\to 
H_i(k+1;R) \to H_{i-1}(k;R)(-1) 
$$
$$
\overset{x_{k+1}}{\to} \cdots \to
H_0(k;R)(-1)\overset{x_{k+1}}{\to} H_0(k;R)\to H_0(k+1;R)\to  0.
$$
The map 
$H_i(k;R)\to H_i(k+1;R)$ is induced by the inclusion of the
corresponding Koszul complexes.
Every homogeneous element $z\in K_i(k+1;R)$
can be uniquely written as $e_{k+1}\wedge z' + z''$
with $z',z''\in K_i(k;R)$. Then
$H_i(k+1;R) \to H_{i-1}(k;R)(-1)$  
is given by sending $[z]$ to $[z']$. 
Furthermore
$H_i(k;R)(-1)\overset{x_{k+1}}{\to} H_i(k;R)$
is the multiplication map with $x_{k+1}$. 
Observe that $H_0(k;R)=R/(x_1,\ldots,x_k) R$.
As noticed above we may assume that the image of $x_1,\ldots,x_n \in S_1$ in $R_1$
is an almost regular sequence for $R$.
In this case the modules $H_i(k;R)$ all have 
finite length for $i>0$. 

We are able to extend Theorem \ref{codim2_cm} to the general case,
which also generalizes the main result in \cite{Go}.

\begin{Thm}
\label{main_2}
Let $I\subset S$ be a graded ideal and $R=S/I$ with $\codim(R) = 2$.
Then 
$$
e(R) \leq (M_1\cdot M_2)/2.
$$
\end{Thm}
\begin{proof}
Let $\xb=x_1,\ldots,x_{n-2}$ and consider
$\Tilde{R}=R/\xb R$. Notice
that by \ref{help_dim} and \ref{help_multi}
we have that
$e(R)\leq e(\Tilde{R})$ and $2=\codim(R)=\codim(\Tilde{R})$.
Observe that $\Tilde{R}=\Tilde{S}/\Tilde{I}$,
where $\Tilde{S}$ is the 2-dimensional polynomial ring 
$S/\xb S$ and $\Tilde{I}=(I+(\xb))/(\xb)$. 
Let
$$
M_i=\max\{j\in \Z\colon \beta^S_{i,j}(R)\neq 0\} \text{ for } i=1,2
$$
and
$$
\Tilde{M}_i=\max\{j\in \Z\colon \beta^{\Tilde{S}}_{i,j}(\Tilde{R})\neq 0\} \text{ for } i=1,2.
$$
We claim that 
\begin{equation}
\label{claim}
\Tilde{M}_1\leq M_1 \text{ and } \Tilde{M}_2\leq M_2.
\end{equation}
Since $\dim(\Tilde{R})=0$, the ring $\Tilde{R}$
is Cohen-Macaulay. Thus
it follows from \ref{codim2_cm} that
$$
e(R)\leq e(\Tilde{R})\leq (\Tilde{M}_1\cdot \Tilde{M}_2)/2 \leq 
(M_1\cdot M_2)/2.
$$
It remains to prove claim (\ref{claim}). 
The first inequality
can easily be seen:
$\Tilde{M}_1$ is the maximal degree of a minimal generator of
$\Tilde{I}$ and
$M_1$ is the maximal degree of a minimal generator of
$I$. Since $\Tilde{I}=(I+(\xb))/(\xb)$ we
get that
$$
\Tilde{M}_1
\leq M_1.
$$
Next we prove the second inequality 
$\Tilde{M}_2\leq M_2$.
Let
$H_{\lpnt}(k;R)$ denote the Koszul homology of $R$
with respect to $x_1,\ldots,x_k \in S$ for $k=1,\ldots,n$
and
$\Tilde{H}_{\lpnt}(l;\Tilde{R})$ denote the Koszul homology of $\Tilde{R}$
with respect to $x_{n-2+1},\ldots,x_{n-2+l} \in \Tilde{S}$ for $l=1,2$.
We denote with
$$
M_{i,k}=\max(\{   j\in \Z \colon  H_{i}(k;R)_j\neq 0  \} \cup \{0\}) \text{ for } i=1,2 \text{ and } k=1,\dots,n.
$$
and 
$$
\Tilde{M}_{i,l}=\max(\{   j\in \Z\colon  \Tilde{H}_{i}(l;\Tilde{R})_j\neq 0  \}\cup \{0\})
\text{ for } i=1,2 \text{ and } l=1,2.
$$
Observe that these numbers are well-defined since all considered modules have finite length.
Note that
$M_{i,n}=M_i$ and $\Tilde{M}_{i,2}=\Tilde{M}_i$ for $i=1,2$.
We have to show that
$$
\Tilde{M}_{2,2}\leq M_{2,n}.
$$
Since $H_0(n-2;R)=\Tilde{R}$
there is the long exact sequence of Koszul homology groups 
$$
\cdots \to 
H_1(n-2;R)
\to 
H_1(n-1;R) 
\to
\Tilde{R} (-1)\overset{x_{n-1}}{\to} 
\Tilde{R} \to \Tilde{R}/(x_{n-1}) \Tilde{R}\to  0.
$$
We also have an exact sequence
$$
0 \to 
\Tilde{H}_1(1;\Tilde{R})
\to
\Tilde{R} (-1)\overset{x_{n-1}}{\to} 
\Tilde{R} \to \Tilde{R}/(x_{n-1}) \Tilde{R}\to  0.
$$
We get a surjective homomorphism
$H_1(n-1;R) \twoheadrightarrow
\Ker(
\Tilde{R} (-1)\overset{x_{n-1}}{\to} 
\Tilde{R})$
and an isomorphism 
$
\Tilde{H}_1(1;\Tilde{R})
\cong
\Ker(
\Tilde{R} (-1)\overset{x_{n-1}}{\to} 
\Tilde{R})$ of graded $K$-vector spaces.
Hence 
$$
 \Tilde{M}_{1,1}\leq M_{1,n-1}.
$$
Next we consider the exact sequence
$$
\cdots \to
H_2(n;R)
\to 
H_1(n-1;R)(-1)
\overset{x_n}{\to}
H_1(n-1;R)
\to 
H_1(n;R) 
\to
\cdots
$$
Since $H_1(n-1;R)_{M_{1,n-1}+1}=0$ 
we have a surjective map 
$$
H_2(n;R)_{M_{1,n-1}+1}
\to 
H_1(n-1;R)_{M_{1,n-1}}.$$
By definition of the number $M_{1,n-1}$ we have that
$H_1(n-1;R)_{M_{1,n-1}}\neq 0$.
It follows that $H_2(n;R)_{M_{1,n-1}+1} \neq 0$
and therefore 
$$
M_{1,n-1}+1 \leq M_{2,n}.
$$

We also have an exact sequence 
$$
0 \to
\Tilde{H}_2(2;\Tilde{R})
\to 
\Tilde{H}_1(1;\Tilde{R})(-1)
\overset{x_n}{\to}
\Tilde{H}_1(1;\Tilde{R})
\to 
\Tilde{H}_1(2;\Tilde{R}) 
\to \cdots
$$
Note that
$\Tilde{H}_1(1;\Tilde{R})_{\Tilde{M}_{1,1}+1}=0$.
Considering the sequence in degree $\Tilde{M}_{1,1}+1$ 
we get an isomorphism 
$\Tilde{H}_2(2;\Tilde{R})_{\Tilde{M}_{1,1}+1} \cong \Tilde{H}_1(1;\Tilde{R})_{\Tilde{M}_{1,1}}\neq 0$
and thus $\Tilde{M}_{1,1}+1 \leq \Tilde{M}_{2,2}$.
In degree $\Tilde{M}_{2,2}$ we obtain
the injective map
$
0\to \Tilde{H}_2(2;\Tilde{R})_{\Tilde{M}_{2,2}}\to \Tilde{H}_1(1;\Tilde{R})_{\Tilde{M}_{2,2}-1}.
$
Since by definition of the number $\Tilde{M}_{2,2}$ we have that 
$\Tilde{H}_2(2;\Tilde{R})_{\Tilde{M}_{2,2}}\neq 0$,
it follows that $\Tilde{H}_1(1;\Tilde{R})_{\Tilde{M}_{2,2}-1}\neq 0$ and therefore
$\Tilde{M}_{2,2} \leq \Tilde{M}_{1,1}+1$.
Hence 
$$
\Tilde{M}_{2,2} =\Tilde{M}_{1,1}+1.
$$
All in all we have shown
that
$$
\Tilde{M}_2=\Tilde{M}_{2,2}=\Tilde{M}_{1,1}+1 \leq M_{1,n-1}+1 
\leq M_{2,n}
=M_2
$$
which is the second part of the desired inequalities of (\ref{claim}).
Thus we proved (\ref{claim}) and
this
concludes the proof.
\end{proof}

%
%
\section{Componentwise linear ideals}  
\label{compsect}
In this section we prove Conjecture \ref{conj2}
for componentwise linear ideals.
We first introduce some notation and recall some definitions. 
(For unexplained notation see \cite{BRHE98}.) 
Given a finitely generated $S$-module $M \neq 0$ 
and $i,j \in \Z$ we denote with
$\beta^S_{i,j}(M)=\dim_K \Tor^S_i(M,K)_j$ the graded Betti numbers of $M$.
Let
$$
\projdim(M)=\max\{i \in \Z \colon \beta_{i,i+j}^S(M) \neq 0 \text{ for some }j\}
$$
be the projective dimension 
and 
$$
\reg(M)=\max\{j \in \Z  \colon \beta_{i,i+j}^S(M) \neq 0 \text{ for some }i\}
$$
be the Castelnuovo-Mumford regularity of $M$. 

For $a=(a_1,\ldots,a_n)\in \N^n$
and a monomial $x_1^{a_1}\cdots x_n^{a_n}\in S$ 
we set $x^a$. 
Let $|a|=a_1+\cdots+a_n$ and $\supp(a)=\{i:a_i \neq 0\}\subseteq [n]=\{1,\ldots,n\}$.
A {\it simplicial complex} $\Delta$
on the vertex set $[n]$
is a collection of subsets of $[n]$
such that $\{i\} \in \Delta$ for $i=1,\ldots,n$, and $F \in \Delta$ whenever 
$F\subseteq G$ for some $G \in \Delta$. 
For $F \in \Delta$ we define $\dim(F)=|F|-1$ where $|F|=|\{i\in F\}|$
and $\dim(\Delta)=\max\{\dim(F)\colon F \in \Delta\}$.
Then $F \in \Delta$
is called an {\it $i$-face} if $i=\dim(F)$. 
Faces of dimension 0, 1 are called {\it vertices}
and {\it edges} respectively.
The maximal faces under inclusion are called the 
{\it facets} of the simplicial complex.
Note that $\emptyset$
is also a face of dimension $-1$. 
For $i=-1,\ldots,\dim(\Delta)$ we define $f_i$ to 
be the number of $i$-dimensional faces of $\Delta$.

We denote with $\Delta^{*}=
\{F : F^{c} \not\in \Delta\}$ the {\it Alexander dual} of $\Delta$ 
where $F^{c}=[n]\setminus F$. This is again a simplicial complex.
For $F=\{i_{1}, \ldots, i_{s}\} \subseteq [n]$ we also 
write $x_F$ for the monomial $\prod_{i\in F}x_i$.
These monomials are also called {\it squarefree monomials}.
Then $K[\Delta]=S/I_{\Delta}$ is the {\it Stanley-Reisner ring} 
of $\Delta$ where 
$$
I_{\Delta}=(x_F: F\subseteq [n],\ F\not\in \Delta)
$$ 
is the {\it Stanley-Reisner} ideal of $\Delta$.
Observe that $\dim(K[\Delta])=\dim(\Delta)+1$. (See \cite{BRHE98} for details.)
At first we relate some of the considered invariants.
For a graded ideal $I \subset S$ let
$$
a(I)=\min\{d \in \Z\colon \beta^S_{0,d}(I)\neq 0 \}
$$
be the {\it initial degree} of $I$.

\begin{Lem}
\label{important}
Let $\Delta$ be a $(d-1)$-dimensional simplicial complex. Then:
\begin{enumerate}
\item
$e(S/I_{\Delta})=\beta^S_{0,a(I_{\Delta^*})}(I_{\Delta^*})$.
\item
$\codim(S/I_{\Delta})=a(I_{\Delta^*})$.
\item
$\projdim(S/I_{\Delta})=\reg(I_{\Delta^*})$.
\end{enumerate}
\end{Lem}
\begin{proof}
Observe that
$F \in \Delta$ is a facet if and only if $x_{F^c}$ is a minimal generator of $I_{\Delta^*}$.
Hence $F$ has maximal dimension $d-1$ if and only if $x_{F^c}$ is a minimal generator 
of $I_{\Delta^*}$ of minimal degree. It follows that
$$
a(I_{\Delta^*})=n-d 
\text{ and } 
\beta^S_{0,a(I_{\Delta^*})}(I_{\Delta^*})=f_{d-1}.
$$

(i):
We know that $e(S/I_{\Delta})=f_{d-1}$. (For example combine 4.1.9 and  5.1.9 in \cite{BRHE98}.)
Thus $e(S/I_{\Delta})=\beta^S_{0,a(I_{\Delta^*})}(I_{\Delta^*})$.

(ii):
This follows from
$$
\codim (S/I_{\Delta})
=
n-\dim (S/I_{\Delta})
=
n-\dim (\Delta) - 1
=
n-d
=
a(I_{\Delta^*}).
$$

(iii): This is a result of Terai in \cite{Te}.
\end{proof}

Recall that an ideal $I \subset S$ is called a {\it monomial ideal} if it is
generated by monomials of $S$. We denote with $G(I)$ 
the unique minimal system of generators for $I$. 
A monomial ideal $I \subset S$ is called {\it squarefree strongly stable},
if it is generated by squarefree monomials such that for all
$x_F\in G(I)$  and $i$ with $x_i|x_F$ we have for all $j<i$ with $x_j\nmid x_F$
that $(x_F/x_i)x_j\in I$.

Note that for a simplicial complex $\Delta$ we have that
$I_{\Delta}$ is squarefree strongly stable if and only if
$I_{\Delta^*}$ is squarefree strongly stable.

We give a new proof for the bound of the multiplicity
in the case of squarefree strongly stable ideals 
which avoids the calculations of the
original proof in \cite{HeSr98}.

\begin{Thm}
\label{first_nice}
Let $\Delta$ be a simplicial complex such that 
$I_{\Delta}$ is a squarefree strongly stable ideal and $c=\codim(S/I_{\Delta})$.
Then 
$$
e(S/I_{\Delta}) \leq (\prod_{i=1}^{c} M_i)/c!.
$$
\end{Thm}
\begin{proof}
Let
$b(S/I_{\Delta})=\max\{i\in \Z \colon \beta^S_{i,i+\reg(S/I_{\Delta})}(S/I_{\Delta})\neq 0\}$.
Since 
$I_{\Delta}$ and $I_{\Delta^*}$ are 
squarefree strongly stable ideals,
it follows from Theorem \ref{cycle_bounded_bettiab} below
that
$$
\beta^S_{i,i+\reg(S/I_{\Delta})}(S/I_{\Delta})\neq 0 \text{ for } i=1,\dots,b(S/I_{\Delta})
\text{ and }
$$
$$
\beta^S_{0,a(I_{\Delta^*})}(I_{\Delta^*})
\leq
\binom{\projdim(I_{\Delta^*})+a(I_{\Delta^*})}{a(I_{\Delta^*})}.
$$
Let 
$p=\projdim(I_{\Delta^*})$. 
We have that
$$
\codim(S/I_{\Delta})
= 
a(I_{\Delta^*})
\leq 
\max\{j \in \Z \colon \beta^S_{p,p+j}(I_{\Delta^*})\neq 0\}
= b(S/I_{\Delta})
$$
where the last equality follows from 
Theorem 2.8 in \cite{BCP}.
(These numbers describe certain ``extremal Betti numbers'' of the considered mo\-dules.)

Hence we get that
$$
M_i=\reg(S/I_{\Delta})+i \text{ for } i=1,\ldots,\codim(S/I_{\Delta}).
$$
Together with the results of Lemma \ref{important} we obtain
$$
e(S/I_{\Delta})
=
\beta^S_{0,a(I_{\Delta^*})}(I_{\Delta^*})
\leq
\binom{\projdim(I_{\Delta^*})+a(I_{\Delta^*})}{a(I_{\Delta^*})}
$$
$$
=
\binom{\reg(S/I_{\Delta})+\codim(S/I_{\Delta})}{\codim(S/I_{\Delta})}
= (\prod_{i=1}^{\codim(S/I_{\Delta})} M_i)/\codim(S/I_{\Delta})!
$$
\end{proof}

For an arbitrary graded ideal
we can prove a weaker bound than the one of Conjecture \ref{conj2}
which was already noticed in
\cite{HeSr98}.
We also get a bound for the codimension of the considered ideal.

\begin{Cor}
\label{thm3}
Let $\chara(K)=0$,
$I \subset S$ be a graded ideal, $R=S/I$ and $c=\codim(R)$. 
Then
\begin{enumerate}
\item
$
c \leq \max\{i \in \Z \colon \beta^S_{i,i+\reg(S/I)}(S/I)\neq 0 \}.$
\item
$e(R) \leq \binom{\reg(R)+c}{c}.$
\end{enumerate}
\end{Cor}
\begin{proof}
Let again
$$
b(S/I)=\max\{ i \in \Z \colon \beta^S_{i,i+\reg(S/I)} (S/I) \neq 0  \}.
$$
By replacing $I$ with the generic initial ideal 
$\Gin(I)$ with respect to the reverse lexicographic order
of $I$ (see for example \cite{Ei95} for details)
we do not change
the multiplicity and the codimension. 
Furthermore by Theorem 2.8 in \cite{BCP} 
also the number $b(S/I)$ does not change. 
This means we may assume that $I$ is a monomial ideal.

By polarization we get a Stanley-Reisner Ideal $I_{\Delta}$
for some complex $\Delta$
with the same Betti diagram as $I$ and also the multiplicity, codimension
do not change.
Hence we may assume that $I=I_{\Delta}$.

Now we replace $I_{\Delta}$ by the Stanley-Reisner ideal of
the associated simplicial complex with respect to symmetric
or algebraic shifting.
Again the multiplicity, codimension and $b(S/I_{\Delta})$
do not change
and we may assume that $I_{\Delta}$ is a squarefree strongly stable ideal.
(See \cite{ArHeHi00} for details on shifting operations.)

In the proof of Theorem \ref{first_nice} 
we showed in fact that for a squarefree strongly stable ideal the desired bounds
of (i) and (ii) hold. 
This concludes the proof. 
\end{proof}

\begin{Rem}\rm
It can also be shown that the bound for the multiplicity of
Corollary \ref{thm3} is valid if $\chara(K)>0$. 
This can be proved analogously to 
the discussion before Corollary 3.8 in \cite{HeSr98}.
\end{Rem}

In a special case we can prove Conjecture \ref{conj2}.
\begin{Cor}
\label{main_result}
Let $I \subset S$ be a graded ideal, $R=S/I$, $c=\codim(R)$ 
and suppose that
$M_i=\reg(R)+i$ for $i=1,\ldots,c$.
Then
$$
e(R) \leq (\prod_{i=1}^{c} M_i)/c!.
$$
\end{Cor}

\begin{Rem}\rm
Corollary \ref{main_result} 
does not imply the upper bound for the multiplicity in Conjecture \ref{conj2} in 
full generality.  
For example even for complete intersections with ideals generated in degree $\geq 2$
the assumptions of the corollary are not satisfied. 

But several known cases besides squarefree strongly stable ideals are included in this result. 
For example the following cases which were originally proved in \cite{HeSr98}
with different proofs for each type of ideal:
\begin{enumerate}
\item
$I$ is a stable ideal.
\item
$I$ is a squarefree stable ideal.
\item
$I$ has a linear resolution.
\end{enumerate}
Next we generalize these results to the case of
componentwise linear ideals.
\end{Rem}

In the following we fix a field $K$ with $\chara(K)=0$.
Recall that
an ideal $I$ is called {\em componentwise linear}, if for all $d\geq 0$ the
ideal $I_{\langle d \rangle}$ has a $d$-linear resolution. 

\begin{Thm}
Let $I \subset S$ be a componentwise linear ideal, $R=S/I$ and $c=\codim(R)$. 
Then 
$$
e(R) \leq (\prod_{i=1}^{c} M_i)/c!.
$$
\end{Thm}
\begin{proof}
Aramova, Herzog and Hibi \cite{ArHeHi} 
proved that
an ideal $I$ 
is componentwise linear if and only if $\beta^S_{i,j}(I)=\beta^S_{i,j}(\Gin(I))$
for all $i,j\in \Z$
where $\Gin(I)$ is the generic initial ideal of $I$ with respect
to the reverse lexicographic order. 
We know that $\Gin(I)$ is stable (see \cite{Ei95}).
Then the Eliahou-Kervaire resolution of $\Gin(I)$
(see also \ref{cycle_bounded_bettiab} below)
and \ref{thm3} (i) 
imply that 
$$
M_i(S/\Gin(I))=\reg(S/\Gin(I))+i \text{ for } i=1,\ldots,\codim(S/\Gin(I)).
$$
Thus we can apply Corollary
\ref{main_result} to conclude the proof.
\end{proof}

We introduce a large class of componentwise linear ideals.
We fix a vector $\ab=(a_1,\ldots,a_n)$ where $2 \leq a_i \leq \infty$.    
The following type of ideal was defined in \cite{GAHIPE} and \cite{Ro01}:
Let $I \subset S$ be a monomial ideal.
$I$ is said to be \it $\ab$-bounded \rm if for
all $x^u \in G(I)$ and all $i \in [n]$
one has $u_i < a_i$. The ideal 
$I$ is called \it $\ab$-stable \rm if, in addition for 
all $x^u \in G(I)$ and all $j \leq m(u)=\max\{i\in [n] \colon u_i\neq 0\}$ with
$u_j<a_j-1$, we have that $x_j x^u/x_{m(u)} \in I$. 
It is easy to see that if $I$ is $\ab$-stable, then 
for all
$x^u \in I$ and all $j \leq m(u)$ with
$u_j<a_j-1$ we have that $x_j x^u/x_{m(u)} \in I$. If $I$ is $\ab$-stable
with $\ab=(2,\ldots,2)$, 
then $I$ is exactly squarefree stable. For
$\ab=(\infty,\ldots,\infty)$ we obtain a stable ideal in the usual sense.

Let $a,b \in \mathbb{Z}$.
We make the convention that $\binom{a}{b}=0$ unless $0\leq b \leq a$.
If $x^u \in S$ with $u \prec \ab$, then we define
$$
l(u)=|\{i \colon u_i=a_i-1,\ i<m(u)  \}|.\label{cycle_bounded_l}
$$ 
The following Theorem was proved in \cite{GAHIPE} and \cite{Ro01}. 
\begin{Thm}
\label{cycle_bounded_bettiab}
Let $I \subset S$ be an $\ab$-stable ideal and $i,j \in \Z$. 
One has, 
independent of the characteristic of $K$,
$$
\beta_{i,i+j}^S(I)=\sum_{x^u \in G(I),\ |u|=j} 
\binom{m(u)-1-l(u)}{i}.
$$
\end{Thm}

As a consequence we are able to determine the regularity 
for $\ab$-stable ideals.
\begin{Cor}
\label{cycle_bounded_abreg}
Let $I \subset S$ be an $\ab$-stable ideal. Then
$$
\reg(I)=\max\{|u| \colon x^u \in G(I)\}.
$$
In particular, if $I$ is generated in degree $d$, then 
$I$ has a $d$-linear resolution.
\end{Cor}

\begin{Cor}
\label{main_result2}
Let $I \subset S$ be an $\ab$-stable ideal, $R=S/I$ and $c=\codim(I)$.
Then 
$$
e(R) \leq (\prod_{i=1}^{c} M_i)/c!.
$$
\end{Cor}
\begin{proof}
Apply \ref{main_result} and \ref{cycle_bounded_bettiab}.
\end{proof}

We can prove a little bit more:

\begin{Thm}
Let $I \subset S$ be an $\ab$-stable ideal.
Then $I$ is componentwise linear.
\end{Thm}
\begin{proof}
For $k\in \N$ let $I_{\leq k}\subset S$ be the ideal which is generated
by all homogeneous polynomials of $I$ of degree at most $k$.

We use the following criterion from \cite{HeReWe}:
a monomial ideal $I$ is componentwise linear if and only if
$\reg(I_{\leq k}) \leq k$ for all $k\in \N$.
Let $I$ be an $\ab$-stable ideal. Then for all $k$ the ideal $I_{\leq k}$ is $\ab$-stable.
By \ref{cycle_bounded_abreg} we have $\reg(I_{\leq k})\leq k$. This concludes the proof.
\end{proof}

\end{document}